%% file: jet_A0.tex
\documentclass[12pt]{article}
\usepackage{amsmath}
\usepackage{amssymb}

\title{\vskip -45pt\bf {}}
\author{ {}}
\date{\ }

\title{\vskip -45pt\bf {In the heart of representable 
metric jets}\footnote{recently published in {\em Diagrammes} (suppléments aux volumes 67+68, 2012, PARIS,pp 33-52)} }
\author{\Large{\bf Elisabeth Burroni}}

\usepackage[a4paper]{typearea}

\usepackage{t1enc}

\input xypic
\usepackage[all]{xy}
\usepackage{picinpar, graphicx}

\input{macro}

\def\tang{\succ\!\!\!\prec}

\begin{document}

\maketitle
\vskip 30pt
\textit{This article is dedicated to Andrée Ehresmann}
\thispagestyle{empty}
\vspace{5mm}

When I first met Andrée, she was 35, and a young dynamic mathematician. She is the one who read the six thesis
(A. and E.Burroni, R.Guitart, Ch.Lair,\break
M. and G.Weidenfeld) defended at Paris 7 in june 1970, Charles Ehresmann being our doctoral advisor.

Since then, I have seen her every week in Paris along with Charles Ehresmann at Ehresmann's Seminar, as well as during the international Category Conferences that she used to organize at the University of Amiens where she was teaching.

I remember the evenings at those Parisian restaurants where we were generously invited by Andrée and Charles Ehresmann along with  visiting Categoricians; that allowed us to meet foreign researchers in the field of Category Theory.

Finally, I am particularly grateful to Andrée for having believed in our joint work (Jacques Penon's and mine) about our metric Differential Calculus.
Here is a synthetic retrospective presentation of this work,
skimming through our previous papers already published in
arXiv \cite{BP1}, TAC \cite{BP2}, the Cahiers \cite{BP3} and JPAA \cite{BP6}.

\section{Introduction}

Here, I aim to immerse myself in the heart of the metric jets, more precisely of those which are representable, restricting myself to the main basic concepts, while going deeper into some notions already mentionned in our previous papers; this will give me the opportunity of lightening the previous texts (including some proofs), while precising some ideas and giving new examples 
(as the bifractal wave function) with a proof at the end of this paper.
Concerning the concrete examples found all along this paper: they play the ``starring role'' in the understanding of 
our metric Differential Calculus!


I give a glossary at the end of this paper 
which briefly recalls some useful notations and definitions
(the first occurrence of a new notion quoted in the text will be 
followed by a *, which suggests to refer to this glossary). 


\section{Metric jets}

So, mainly, as its name shows it, our metric Differential Calculus generalizes the classical Differential Calculus in a context which is {\it a priori} only metric.
In this context, the metric jets play the part of the differential maps.
The proofs of the assertions of this section can be found in 
the first chapter of \cite{BP1}, in \cite{BP2} and in \cite{BP6};
except for the proof of the fact that the metric jet $\cali$ is a good jet (see examples 2.1 below) that can be found at the end of this paper.
In this section, I recall the concepts which are useful 
for the understanding of the next section 3.

\subsection{The category $\J$et}

Let $(M,a)$ and $(M',a')$ be two pointed metric spaces (the chosen points being always assumed to be non isolated).
A \textit{metric jet} (in short \textit{jet}) $\varphi:(M,a)\lra (M',a')$ is an equivalence class (of maps $f:M\lra M'$ which are $LL_a$* and verify $f(a)=a'$) for the equivalence relation (of \textit{tangency} at $a$):
$f\tang_a g$ if $f(a)=g(a)$ and
$\lim_{x\rightarrow a}\frac{d(f(x),g(x))}{d(x,a)}=0$.
Indeed, when $M$ and $M'$ are n.v.s.*, we come accross 
the usual notion of tangency at a point: more precisely, if $f$ is 
$Dif\! f_a$*, we have $f\tang_a$ A$f_a$ where 
A$f_a(x)=f(a)+$d$f_a(x-a)$ is the continuous affine map which is tangent to $f$ at $a$.

For such a jet $\varphi:(M,a)\lra(M',a')$, we are interested in 
its \textit{lipschitzian ratio}
$\rho(\varphi)=\inf\{k>0\ |\ \exists f\in\varphi,\ f\quad  
k$-$LL_a\}$.
This ratio verifies the inequality
$\rho(\varphi'.\varphi)\leq\rho(\varphi')\rho(\varphi)$,
the composition of jets being defined just below.
  
The local lipschitzianity is a sufficient condition for composing the jets; indeed, if 
$\xymatrix{
M\ar@<1ex>[rr]^{f}\ar@<-1ex>[rr]_{g}
&&M'\ar@<1ex>[rr]^{f'}\ar@<-1ex>[rr]_{g'}&&
M''
}$,
where $M,M',M''$ are metric spaces respectively pointed by $a,a',a''$, the quoted maps verifying $f(a)=g(a)=a'$, 
$f'(a')=g'(a')=a''$, and $f,g$ (resp. $f',g'$) being $LL_a$   
(resp. $LL_{a'}$),
then, we have the implication:
$f\tang_{a} g$ and 
$f'\tang_{a'} g'\Lra f'.f\tang_a g'.g$. 
So, the jets are the morphisms of a category $\J$et whose objects are the pointed metric spaces ($\varphi'.\varphi$ is the jet containing $f'.f$ where $f\in\varphi$ and $f'\in\varphi'$).
This category $\J$et is cartesian and enriched in  
$\M$et*, a ``well-chosen'' category of metric spaces
(whose morphisms are $LSL$* maps).
Thus, we can speak of the distance $d(\varphi,\psi)$
when $\varphi,\psi\in\J$et$((M,a),(M',a'))$ as being the ``quasi-distance'' $d(f,g)$* where $f\in\varphi$ and $g\in\psi$.  

More precisely, if we denote
$\L$L$((M,a),(M',a'))$ the set of maps $f:M\lra M'$\break
which are $LL_a$ and verify $f(a)=a'$ (providing ``Hom''
for a cartesian category denoted $\L$L),
and $q$ the canonical surjection $\L$L$((M,a),(M',a'))\lra\J$et$((M,a),(M',a'))\break
=\L$L$((M,a),(M',a'))/\tang_a$,
then the quasi-distance* on
$\L$L$((M,a),(M',a'))$
factorizes through the quotient
giving a distance, 
defined by $d(q(f),q(g))=d(f,g)$.

We notice that, for all $\varphi\in \J$et$((M,a),(M',a'))$, we have
$d(\varphi,\calo_{aa'})\leq\rho(\varphi)$, where $\calo_{aa'}
:(M,a)\lra(M',a')$ is the jet containing the constant map on $a'$.
A jet $\varphi$ is said to be a \textit{good jet} if the previous inequality is an equality.

\begin{exams}
{}
\end{exams}

Here, except for example 5 (for which $M=]-1,1[$ and $M'=\R$), we have $M=M'=\R$; and everywhere $a=a'=0$. We set $\calo=\calo_{00}$.

\qquad 1) $\calo$ is the jet of all the $LL_0$ maps which are tangent at 0 to the constant 
map on 0; it is a good jet with $\rho(\calo)=0$.

\qquad 2) $\calv$ is the jet containing the absolute value 
$v(x)=|x|$; 
\vspace{6mm}

{\centerline
   {\includegraphics[width=8cm]{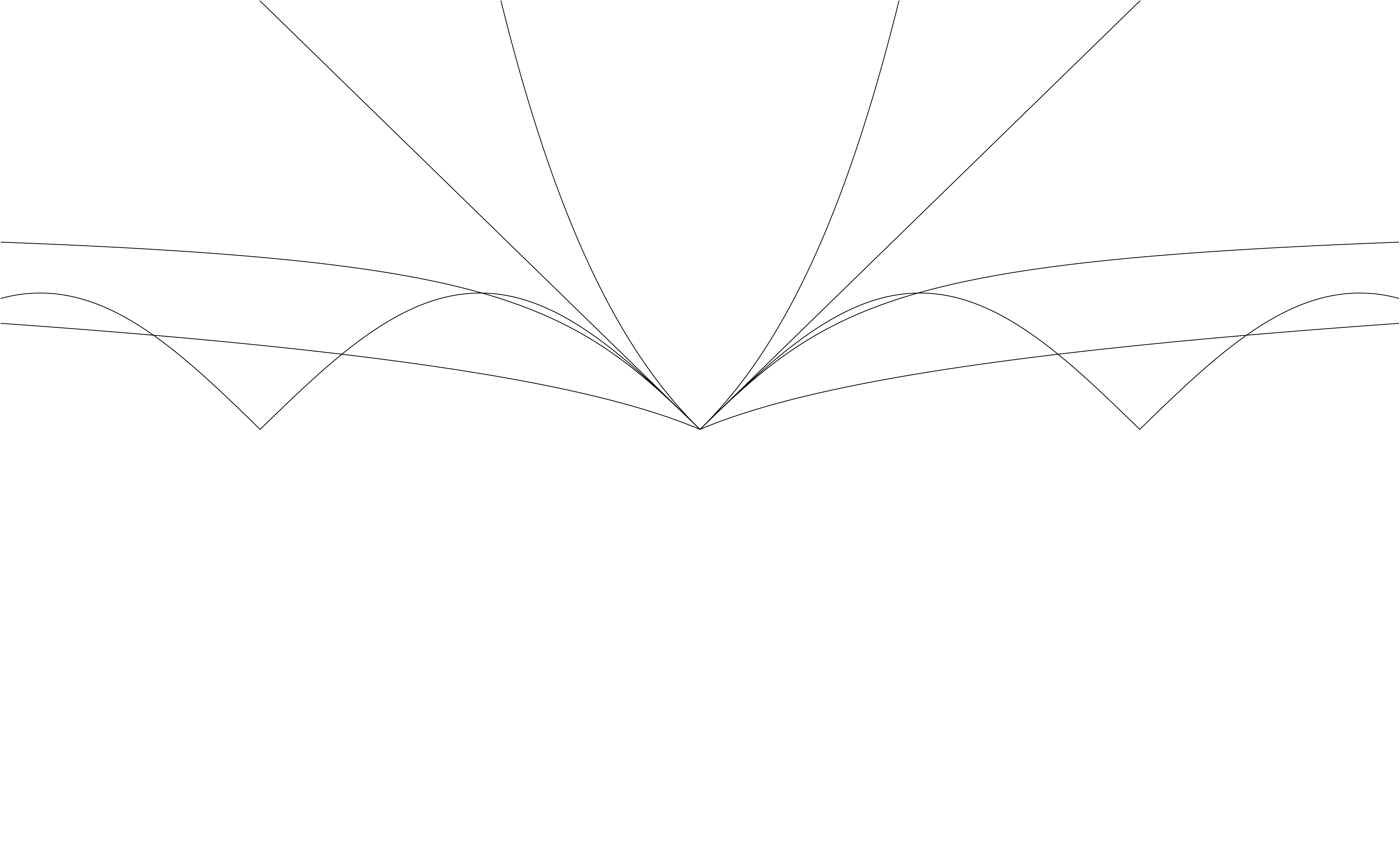}}
  }
  
\vspace{-18mm}  
This jet $\calv$ contains the functions (considering Taylor expansions at order 1) $\exp|x|-1$, $|x|$, $\sin|x|$, $\log(1+|x|)$, $Arctg |x|$\ \dots\ \textit{etc}.
It is a good jet since
$d(\calv,\calo)=\break
d(v,0)$*
$=\lim_{r\rightarrow 0}\sup
\{\frac{v(x)}{|x|}\ |\ 0\not=|x|\leq r\}=1\leq\rho(\calv)\leq 1$, this last inequality being due to the fact that $v$ is 1-lipschitzian.

3) $\calg$ is the jet of the Giseh function
$g(x)=d(x,K_\infty)$ where $K_\infty=\bigcup_{n\in\N}3^n\K$
and $\K$ is the triadic Cantor set. This jet is a good jet with $\rho(\calg)=1$.

4) $\calf$ is the jet of the fractal wave function 
$\xi(x)=x\sin\log|x|$ if $x\not=0$, $\xi(0)=0$; this jet is not a good one since $1=d(\calf,\calo)<\rho(\calf)=\sqrt 2$.  

5) $\cali$ is the jet of the uncanny function (in french ``insolite'')
$Ins:]-1,1[\lra\R:x\mapsto x\sin\log|\log|x||$ if $x\not=0$, $Ins(0)=0$;
we prove (see Proof 1 at the end of this paper) that this jet $\cali$ is a good one,
with $\rho(\cali)=1$.

We will meet again these examples farther;
in particular, one can
 find the graphs of $g$ and $\xi$ in examples 3.2).
 
\vspace{3mm}
We conclude these examples with the jet $\calj_a$ containing the canonical injection $j:V\hookrightarrow M$, where $V$ is a neighborhood  of $a$ in a metric space $M$;
this jet $\calj_a:(V,a)\lra (M,a)$ is an isomorphism in $\J$et,
its inverse $\calj_a^{-1}$ being the jet of
the map $s:M\lra V$ defined by $s|_V=id_V$ and $s(x)=a$ when 
$x\notin V$. These jets are good jets with
$\rho(\calj_a)=\rho(\calj_a^{-1})=1$.

Let $(M,a)$ and $(M',a')$ be pointed metric spaces,
$V$ (resp. $V'$) a neighborhood of $a$ in $M$ 
(resp. of $a'$ in $M'$).
If $\varphi\in\J$et$((V,a),(V',a'))$, we denote 
$\Gamma(\varphi)$ the following composite jet (that we call the
\textit{stretching} of $\varphi$ to $M$):

$$\xymatrix{
(M,a)\,
\ar[r]^{\!\!\!
\calj_a^{-1}}& (V,a)\,
\ar[r]^{\!\!\!\!\!\!\!\varphi}&\,(V',a')\,
\ar[r]^{\,\,\,\,\calj_{a'}}&\,(M',a')
}$$

This defines an isometry $\Gamma:\J$et$((V,a),(V',a'))\lra
\J$et$((M,a),(M',a'))$ which verifies
$\rho(\Gamma(\varphi))=\rho(\varphi)$
(since $\rho(\Gamma(\varphi))\leq\rho(\calj_{a'})
\rho(\varphi)\rho(\calj_a^{-1})=\rho(\varphi)$;
same for the inverse inequality), so that
$\Gamma(\varphi)$ is a good jet iff $\varphi$ is a good jet. 

\subsection
{Tangent jets}

Untill now, we only have spoken
of jets for maps which are locally lipschitzian at a point.
More generally, we can associate a jet to a map which is tangentiable
at a point; the tangentiable maps are natural generalisations of the
differentiable maps.

Let $M,M'$ be two metric spaces, $f:M\lra M'$ a map and $a\in M$.
We say that $f$ is \textit{tangentiable} at $a$ (in short $Tang_a$) if there exists
an $LL_a$ map  $g:M\lra M'$ such that $f\tang_a g$.
If $f$ is $Tang_a$, the set 
$\{g:M\lra M'\, |\, g\ $is $\ LL_a\ $and$\ f\tang_a g\}$ is a jet $(M,a)\lra(M',f(a))$ which is denoted T$f_a$ and called the \textit{tangent jet} of $f$ at $a$.
By definition, 
$f\ LL_a\Llra f\ Tang_a$ with
$f\in$ T$f_a$ (i.e T$f_a=q(f)$); in fact,\break
$f$ $LL_a\Lra f$ $Tang_a\Lra
f$ $LSL_a\Lra f\ C^0_a$*.

We have a composition of the tangent jets for composable tangentiable maps:
T$(g.f)_a=$T$g_{f(a)}.$T$f_a$ if $f$ is $Tang_a$ and $g$ is $Tang_{f(a)}$.

\begin{exams}
{}
\end{exams}

1) $f$ $Dif\! f_a$* $\Lra$ $f$ $Tang_a$,  
where T$f_a$ is the jet containing the continuous affine map A$f_a$ 
tangent to $f$ at $a$.
 
2) All the examples 2.1 are lipschitzian but not $Dif\! f_0$: we have\break
$\ v\in\calv=$T$v_0$,  
$\ g\in\calg=$T$g_0$,
$\ \xi\in\calf=$T$\xi_0$,
$\ Ins\in\cali=$T$Ins_0$. We also have
$\ j\in\calj_a=$T$j_a\ $ \dots $\ $ \textit{etc}.

3) $f_1(x)=x\sin\frac{1}{x}$ if $x\not=0$, $f_1(0)=0$, is 
not $Tang_0$, although obviously $LSL_0$.

4) $f_2(x)=x^2\sin\frac{1}{x^2}$ if $x\not=0$, $f_2(0)=0$, 
is $Tang_0$ (since it is $Dif\! f_0$) but not $LL_0$
(since $\lim_{k\rightarrow+\infty}f'_2(\frac{1}{\sqrt{2k\pi}})=-\infty$);
so that the tangent jet T$(f_2)_0$ exists, but
$f_2\notin$T$(f_2)_0$.

\subsection{Inside n.v.s.*}

Untill now we have contented ourselves with a purely metric context; from now on, we will consider all the previous new notions
in the n.v.s. frame (which will provide new concepts for such a classical frame). We denote $E$, $E'$ \dots \, such n.v.s.. 

First, we notice that, like $\L$L$((E,a),(E',0))$,  the set
$\J$et$((E,a),(E',0))$ is a vector space
(since the vector space $(E',0)$ is also a
vector space internally in $\J$et: see examples 2.3 below);
and the canonical surjection
$q:\L$L$((E,a),(E',0))\lra\J$et$((E,a),(E',0))$
is a linear map.
In fact, $\J$et$((E,a),(E',0))$ is a n.v.s., its distance deriving from a norm $\|\varphi\|=d(\varphi,\calo_{a0})$.
Thus $\varphi:(E,a)\lra(E',0)$ is a good jet iff $\rho(\varphi)=\|\varphi\|$.

Notably, we find the following good jets:
if $l:E\lra E'$ is a continuous linear map, then $l$ is $Tang_0$ with $l\in$T$l_0:(E,0)\lra(E',0)$, and
the restriction  
$\L(E,E')\hookrightarrow\L$L$((E,0),(E',0))
\buildrel q\over\lra\J$et$((E,0),(E',0)):\ l\ \mapsto$ T$l_0$ 
is a linear isometric embedding
($\L(E,E')$ being the set of all continuous linear maps $E\lra E'$, equipped with the operator norm
$\|l\|_{op}=\sup_{x\not=0}\frac{\|l(x\|}{\|x\|}$). So we have
$\|l\|_{op}=\|q(l)\|=\|$T$l_0\|=d($T$l_0,\calo)\leq \rho($T$l_0)\leq
\|l\|_{op}$, the last inequality being due to the well-known fact that $l$ is $\|l\|_{op}$-lipschitzian; this implies that T$l_0$
is a good jet.

\begin{exams}
{}
\end{exams}

1) If $E$ is a n.v.s., we denote
$\sigma:E\times E\lra E:(x,y)\mapsto x+y$ and 
$m_\lambda:E\lra E:x\mapsto\lambda x$ the continuous linear operations of $E$;
Then, T$\sigma_{(0,0)}:(E,0)^2\lra(E,0)$ and T$m_\lambda:(E,0)\lra(E,0)$ are good jets,
respectively denoted + and $\mu_\lambda$.
The data of these two jets confers on $(E,0)$ a structure of vector space, internally in $\J$et.

2) The translation $\theta_{ab}:E\lra E:x\mapsto x+b-a$ provides a jet T$(\theta_{ab})_a:(E,a)\lra(E,b)$ denoted $\gamma_{ab}$ which verifies $\rho(\gamma_{ab})\leq 1$ and is 
invertible in $\J$et with $\gamma_{ab}^{-1}=\gamma_{ba}$. 
If $\varphi\in\J$et$((E,a),(E',a'))$, we denote
$\Omega(\varphi)$ the following composite jet (that we call
the \textit{translate} of $\varphi$ in 0):

$$\xymatrix{
(E,0)\,
\ar[r]^{\!\!\!
\gamma_{0a}}& (E,a)\,
\ar[r]^{\!\!\!\!\!\!\!\varphi}&\,(E',a')\,
\ar[r]^{\!\!\!\,\gamma_{a'0}}&\,(E',0)
}$$

This defines an isometry
$\Omega:\J$et$((E,a),(E',a'))\lra\J$et$((E,0),(E',0))$ which verifies 
$\rho(\Omega(\varphi))=\rho(\varphi)$, so that
$\Omega(\varphi)$ is a good jet iff $\varphi$ is a good jet.

\subsection{Tangentials}

\vspace{2mm}
If a map $f:U\lra U'$ is $Tang_a$, where $U$ and $U'$ are open subsets of $E$ and $E'$ respectively, $a\in U$,  
we denote t$f_a:(E,0)\lra(E',0)$ the following composite jet (that we call the \textit{tangential} of $f$ at $a$):

$\xymatrix{
(E,0)\,\ar[r]^{\!\!\!\!\gamma_{0a}} &\,(E,a)\,
\,\ar[r]^{\!\!\! \calj_a^{-1}}&\, (U,a)\,
\ar[r]^{\!\!\!\!\!\!\!\textup{T}f_a}&\,(U',f(a))
\ar[r]^{\!\!\!\!\!\calj_{f(a)}}&\,(E',f(a))\,
\ar[r]^{\,\,\,\,\gamma_{f(a)0}}&\,(E',0)
}$

In other words, this jet t$f_a=\Omega(\Gamma($T$f_a))$ is a ``streched translate at $0$'' of the tangent jet T$f_a$ in $\J$et.
Now, if $f$ is tangentiable at every point of $U$,
it provides a map t$f:U\lra\J$et$((E,0),(E',0)):\, x\, \mapsto$
t$f_x$, which is called the \textit{tangential} of $f$.

\subsection{linear jets}

I complete this section with the notion of linear jet.

A jet $\varphi:(E,0)\lra(E',0)$ is said to be \textit{linear}
if  the following two
diagrams commute in the category $\J$et:

$$\xymatrix{
(E,0)^2\ar[d]_{+}\ar[r]^<<<<<{\varphi^2}&(E',0)^2
\ar[d]^{+}&&
(E,0)\ar[d]_{\mu_\lambda}\ar[r]^<<<<<{\varphi}&(E',0)
\ar[d]^{\mu_\lambda}
\\
(E,0)\ar[r]_<<<<<{\varphi}&(E',0)&&
(E,0)\ar[r]_<<<<<{\varphi}&(E',0)
}$$
\vspace{1mm}
\noindent where the jets 
+ and $\mu_\lambda$ have been defined in
examples 2.3.

These two commutative diagrams merely mean that the jet  $\varphi$ is linear internally in $\J$et.
The set of the linear jets $(E,0)\lra(E',0)$ is denoted $\Lambda(E,E')$; 

If $l:E\lra E'$ is a continuous linear map, its tangent jet
T$l_0$ is a good linear jet (just apply the composition of the tangent jets to the equalities
$l.\sigma=\sigma.l^2$ and $l.m_\lambda=m_\lambda.l$); ex: + and $\mu_\lambda$ are good linear jets, so that the set $\Lambda(E,E')$ is a sub-n.v.s. of $\J$et$((E,0),(E',0))$;
more precisely,
the linear  isometric embedding $\L(E,E')\hookrightarrow 
\J$et$((E,0),(E',0)):l\mapsto$ T$l_0$ factorizes through
$\Lambda(E,E')$.

\begin{exams}
{}
\end{exams}

1) As previously said, the continuous linear maps give rise to linear jets; but, as the following example shows it, a linear jet is not necessarily the jet of a continuous  linear map
(however, we will see in theorem 3.11 that it can be true in a specific context).

2) The tangent jet $\cali=$T$Ins_0$ (defined in examples 2.2) is not linear ($]\! -\! 1,1[$ being not a vector space); but its streching jet
$\Gamma(\cali)$ to $\R$:
$\!\!\xymatrix{
(\R,0)
\ar[r]^{\!\!\!\!\!
\calj_0^{-1}}& (]\! -\! 1,1[,0)
\ar[r]^{\ \ \textup{T}Ins_0}&(\R,0)
}$ is a linear jet.
In fact, $\Gamma(\cali)=$T$\overline{Ins}_0$ where $\overline{Ins}$ is the extent of $Ins$ to $\R$ (giving the value 0 on $]-1,1[^c$); this extent keeps all the local properties of $Ins$ at 0.

Let us stand still for a while on this uncanny function
$Ins$,
just to have a better understanding of the notion of linear jet.

The commutativity of the two square diagrams expressing the linearity of  the jet $\Gamma(\cali)$ simply means that
$Ins.\sigma\tang_{(0,0)} \sigma. Ins^2$ and
$Ins. m_\lambda\tang_0 m_\lambda. Ins$; \textit{i.e} that
the uncanny function $Ins$ verifies:

\vspace{-5mm}
$$\lim_{(x,y)\rightarrow (0,0)}
\frac{Ins(x+y)-Ins(x)-Ins(y)}{\|(x,y)\|}=0\qquad \\
\qquad \lim_{x\rightarrow 0}
\frac{Ins(\lambda x)-\lambda Ins(x)}{x}=0$$

This could be expressed saying that $Ins$ is  ``linear at the      limit'' at 0; let us have a look on the graph of $Ins$, just to
get a good idea of this ``limit linearity''.
If, at first sight, it may seem rather
simple, the appearances are however misleading:

\vspace{2mm}

\centerline{
\includegraphics[width=4cm]{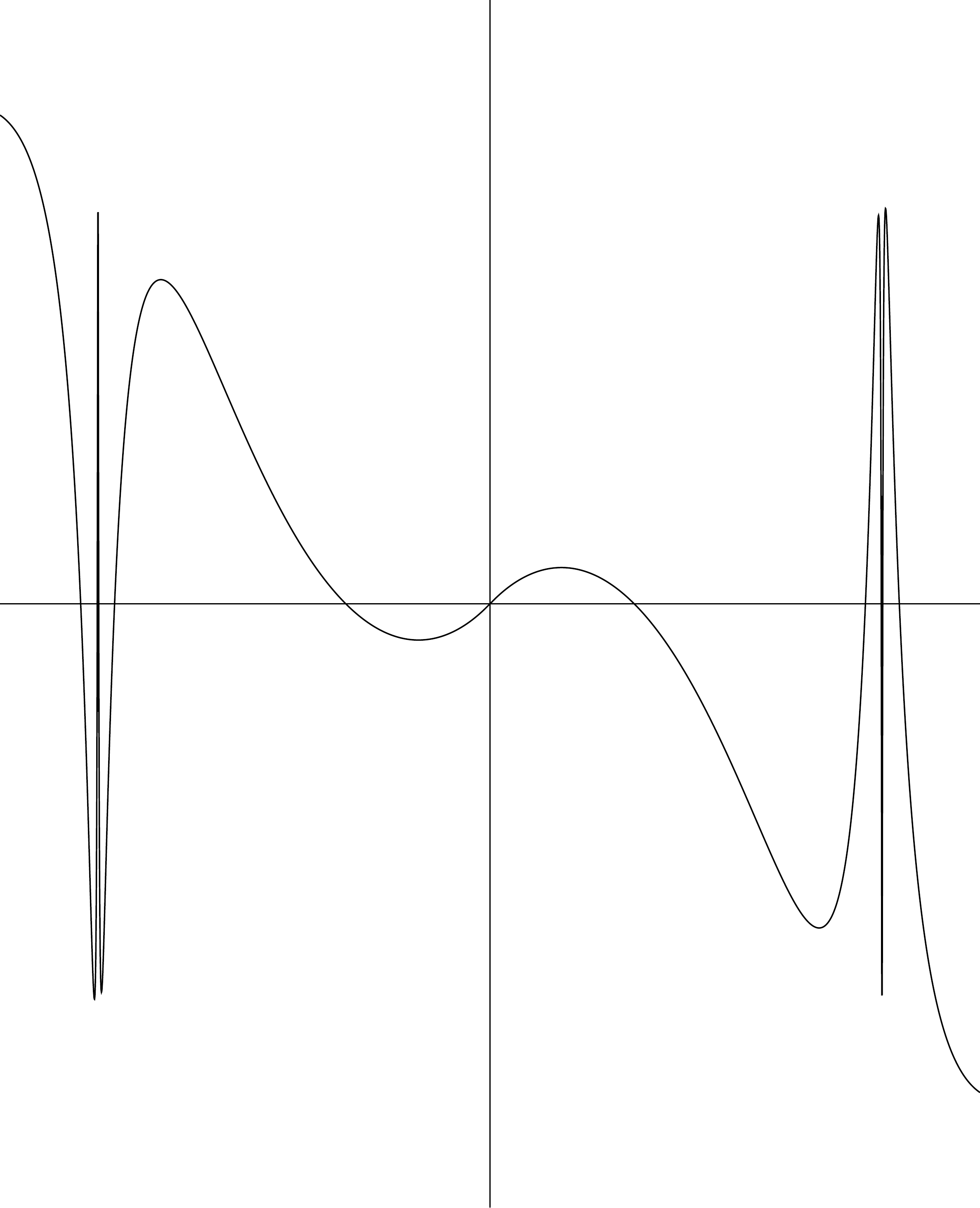}\qquad\
\includegraphics[width=25mm]{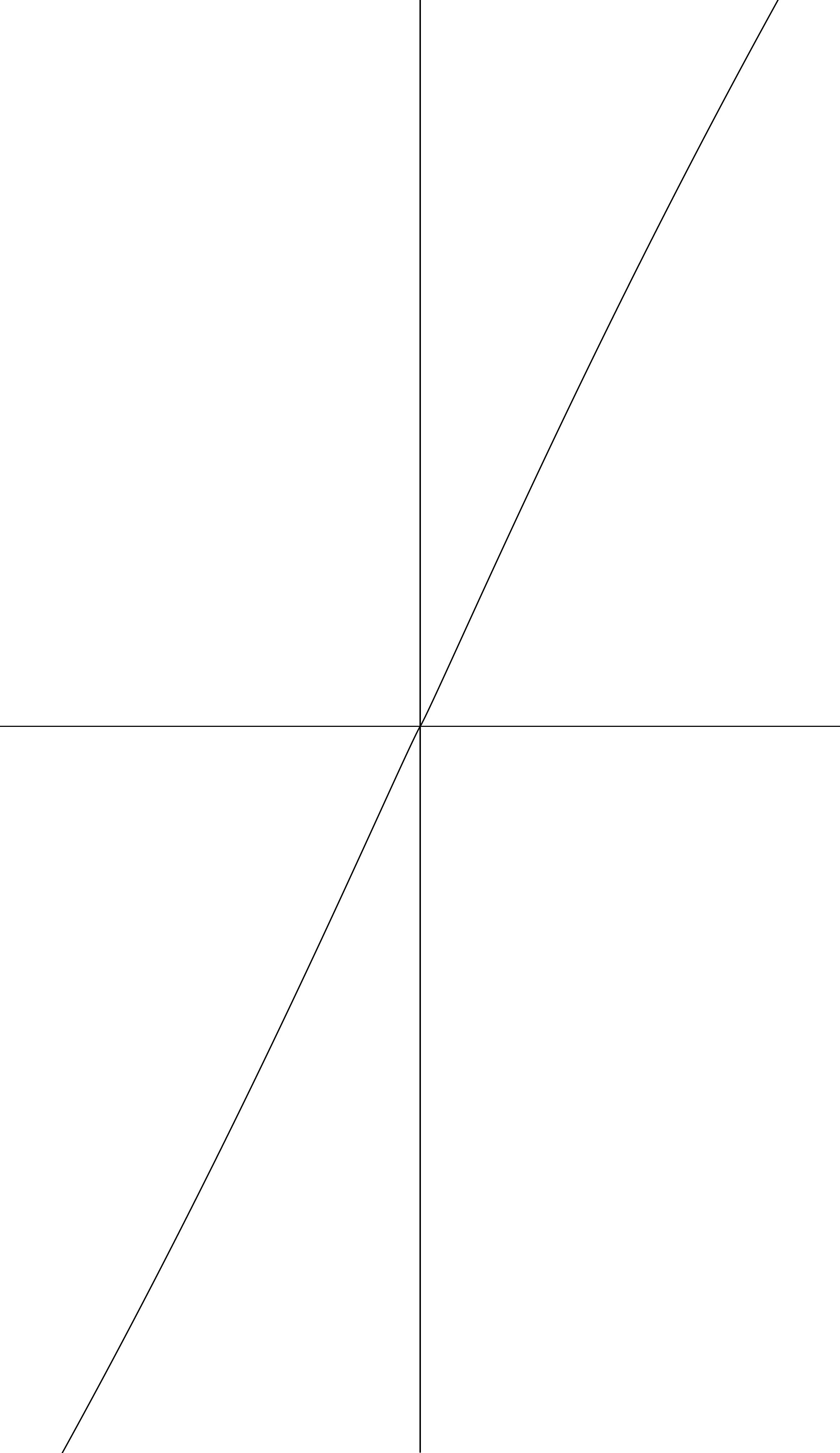}
\includegraphics[width=25mm]{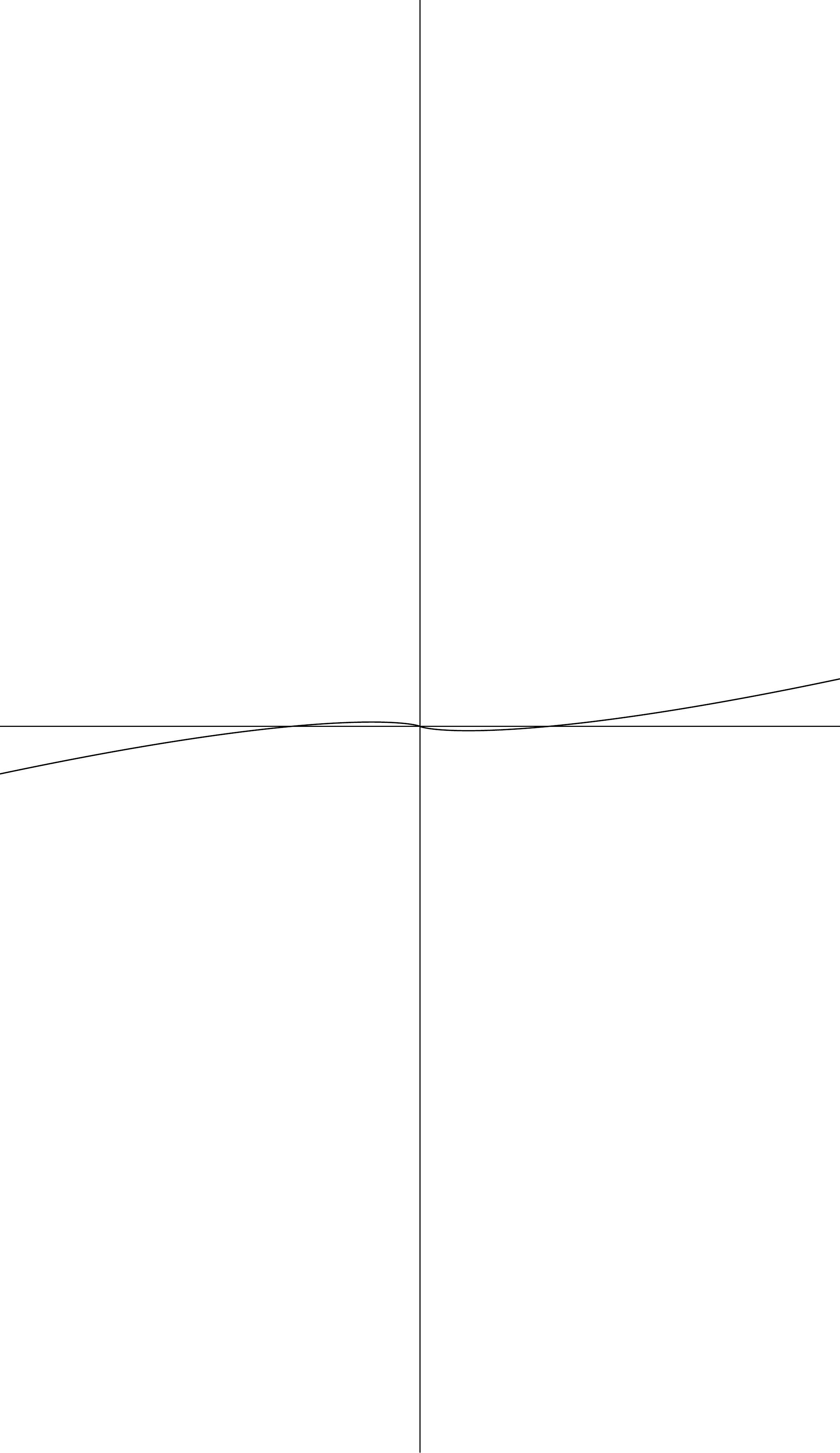}
\includegraphics[width=25mm]{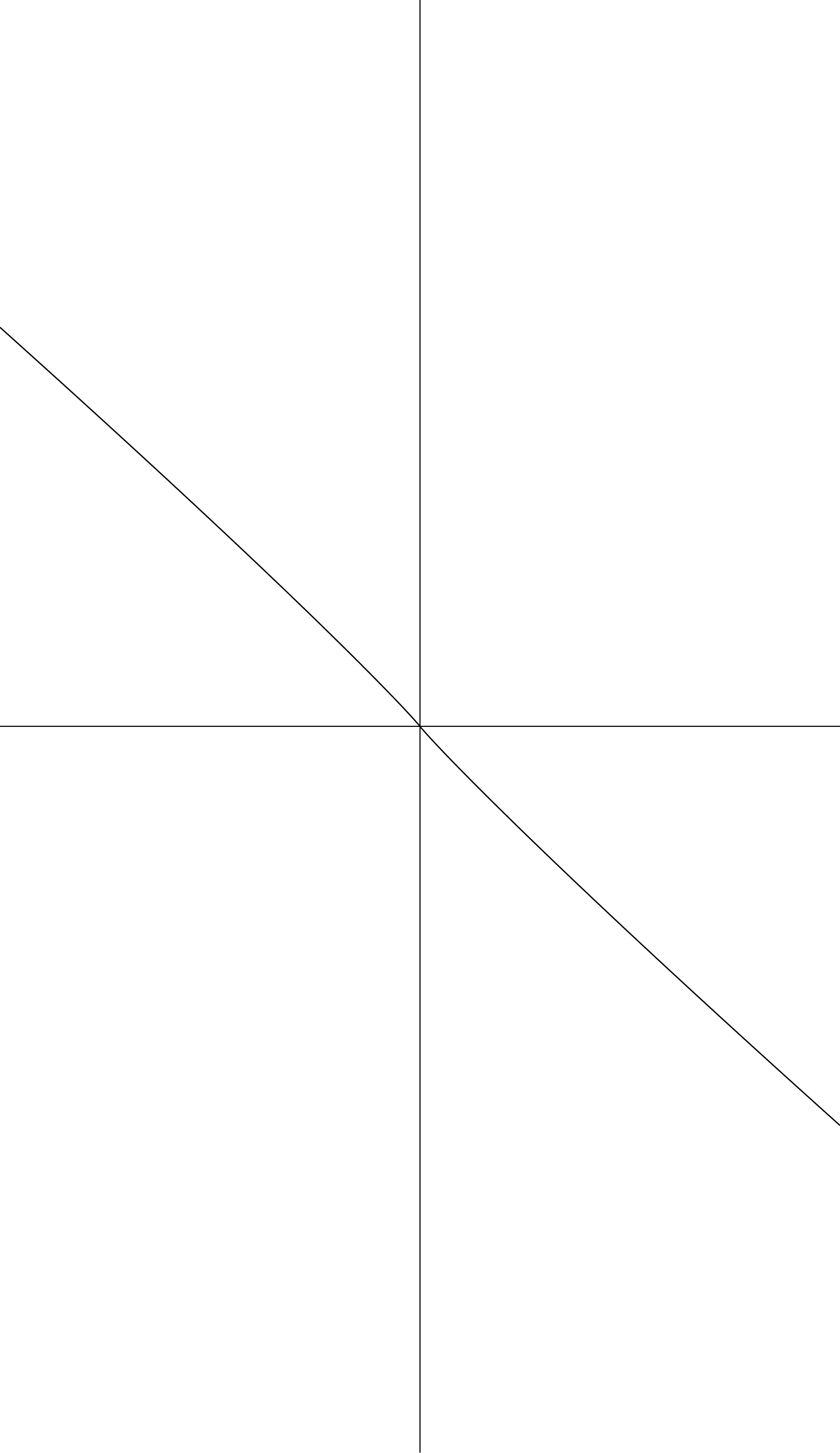}
}

\noindent considering
more and more  powerful zooms on 0, 
we notice that, getting closer and closer to 0,
the slope is constantly changing (which only means that 
$Ins$ is not $Dif\! f_0$);
the most important thing being that more we get closer to 0,   more the function ``is''  rectilinear.
And, indeed, ``more and more rectilinear'' 
can be expressed saying ``linear at the          limit''!

\subsection{Tangentially linear maps}

The notion of linear jet gives us the occasion of defining a new generalization of differentiable maps (of course still in the n.v.s. context).

If $U$ and $U'$ are open subsets of $E$ and $E'$ respectively, and if $a\in U$, a map $f:U\lra U'$
is said to be \textit{tangentially linear} at $a$ (in short $TL_a$) if $f$ is $Tang_a$ and if its tangential at $a$ t$f_a:(E,0)\lra(E',0)$
is linear
(we could say that the tangent jet T$f_a$ is affine).

Now, if $f$ is tangentially linear at every point of $U$,
its tangential t$f:U\lra\J$et$((E,0),(E',0)):x\mapsto$ t$f_x$ factorizes through $\Lambda(E,E')$;
we denote $\lambda f$ the restriction of t$f$ to
$U\lra\Lambda(E,E')$.

\begin{exams}
{}
\end{exams}

1) $f\ Dif\! f_a\ \Lra\ f\ TL_a$ with
d$f_a\in\ $t$f_a=$T(d$f_a)_0$.
For example, $f_2(x)=x^2\sin\frac{1}{x^2}$ if $x\not=0$, $f_2(0)=0$,
is $TL_0$ with t$(f_2)_0=$T$(f_2)_0=$T(d$f_2)_0=0$.

2) The uncanny function $Ins:]-1,1[\lra\R$ is $TL_0$,
but not $Dif\! f_0$.

\vspace{3mm}
We end this section with the following 
$\mathbf{local\ inversion\ theorem}$:

\begin{theo}

Let $f:U\lra U'$ be a $CTL$* map, where $U$ and $U'$ are open subsets of $E$ and $E'$ (supposed to be Banach spaces)
respectively, and $a\in U$.
We assume that there exists an invertible germ $G:E\lra E'$ in $\G{\rm{CTL}}$*
such that $G\subset\lambda f_a$. Then, there exists an open neigborhood $V$ of $a$ in $U$ such that $f(V)$ is open in $E'$ and that the
restriction of $f$ to $V\lra f(V)$ is invertible in $\C$\rm{TL}*.
\end{theo}

Here, the invertible germ $G$ verifying $G\subset\lambda f_a$ plays the part of the invertible differential
d$f_a$ ($f$ being then supposed to be of class $C^1$)
of the classical local inversion theorem.

\section{Representable metric jets}

I now come to the heart of my subject, still in the simplifying 
n.v.s. framework (even though it is possible to work in more general specific metric spaces that we call $\Sigma$-contracting spaces*).
The proofs of the assertions of this section can be found in the second chapter of \cite{BP1} and in \cite{BP3}; except for 
the proofs of the theorem 3.11 and of the fact that the bifractal wave function is not neofractal at 0 (see examples 3.14) that can be found at the end of this paper.
All the basic concepts used here have been recalled in the previous section 2.

\subsection {Valued monoid}

A monoid $\Sigma$ (whose law is denoted multiplicatively) is 
called a {\textit {valued monoid}} if it is equipped with a 
specific element 0 and with a homomorphism
$v:\Sigma\lra\R_+$ verifying the two conditions:
$v(t)=0\Llra t=0$ and
$\exists t\in\Sigma\quad (0<v(t)<1)$.\break
A map $\sigma:\Sigma\lra\Sigma'$ is said to be a morphism of valued monoids
if it verifies $\sigma(0)=0$ and $v(t)=v'(\sigma(t))$ for all $t\in\Sigma$. 

\vfill\eject

\begin{exams}
{}
\end{exams}

We will consider the two following examples
of morphisms of valued monoids:\break 
$\N'_k\hookrightarrow\R_+\hookrightarrow\R$ 
(with
$0<k<1$), where $\N'_k=\{k^n\,|\, n\in\N\}\cup\{0\}$ is
valued by $id_{\N'_k}$; $\R_+$  valued
by $id_{\R_+}$, and $\R$ valued by $v(t)=|t|$.

\vspace{3mm}
If $E$ is 
a n.v.s. (and $a\in E$), 
we denote it $E_a$ if we consider it as being centred in $a$: 
this $E_a$ is a $\R$-vector space with $0_a=a$, $x+_ay=a+((x-a)+(y-a))$
and $\lambda._a x=a+\lambda(x-a)$; it is even a n.v.s., setting
$\|x\|_a=\|x-a\|$. Of course, $E_0=E$ with its given norm.

Now, every valued monoid $\Sigma$ provides $E_a$ with a 
new canonical external operation by setting $t\star_a x=a+v(t)(x-a)$;
which, besides the usual properties of external operations, verifies 
$0\star_a x=a$ for all $x\in E_a$ and
$t\star_a a=a$ for all $t\in\Sigma$;
and is compatible with the norm on $E_a$
\textit{i.e} verifies, for all $t\in\Sigma$ and $x\in E_a$, $\|t\star_a x\|_a=v(t)\|x\|_a$,
which implies $\lim_{n\rightarrow +\infty}t^n\star_a x=a$ when $t\in\Sigma$ verifies $0<v(t)<1$.
Besides, for every morphism of valued monoids  $\sigma:\Sigma\lra\Sigma'$,
we have $t\star_a x=\sigma(t)\star'_a x$ for 
all $t\in\Sigma$ and $x\in E_a$.

\subsection{Homogeneous maps}

In all that follows, $\Sigma$ is a valued monoid.

A map $h:E_a\lra E'_{a'}$ is said to be 
$\Sigma$-\textit{homogeneous} if it verifies
$h(t\star_a x)=\break
t\star'_{a'}h(x)$ for all $t\in\Sigma$ and $x\in E_a$; such an homogeneous map verifies\break
$h(a)=h(0\star_a a)=0\star'_{a'}h(a)=a'$.
Thanks to the morphisms of valued monoids $\Sigma\buildrel v\over\lra \R_+\hookrightarrow\R$, we have
the implications:
$\R$-homogeneous $\Lra$ $\R_+$-homogeneous $\Lra$
$\Sigma$-homogeneous for all $\Sigma$.
 
If we consider n.v.s. centred in 0, then
$h:E\lra E'$ is $\R_+$-homogeneous if it 
verifies $h(tx)=th(x)$ for all $t\in\R_+$ and $x\in E$, 
\textit{i.e} $h$  is 
positively-1-homogeneous; and
$h:E\lra E'$ is $\N'_k$-homogeneous if 
it verifies the fractal property $h(kx)=kh(x)$ for 
all $x\in E$; it is why we say $k$-\textit{fractal} instead of 
$\N'_k$-homogeneous.

Why fractal? Because of the equivalence:
$(x,y)\in Graph(h)$ iff $(kx,ky)\in Graph(h)$, meaning that 
$Graph(h)$ remains identical to itself when we zoom into 
0 with a ratio $k$ (this process being iterated 
for an infinity of times).
We can have an approximative idea of a
fractal function $h:\R\lra\R$, by considering 0 as a point at the infinity (i.e at the unreachable horizon point),
the graph of $h$, being then seen in perspective, infinitely decreasing towards this horizon point, and still remaining itself, but thinner and thinner.

\vfill\eject

\begin{exams}
{}
\end{exams}

1) Linear $\Lra$ $\R_+$-homogeneous $\Lra$ $\Sigma$-homogeneous, for all $\Sigma$.

2) The function $v(x)=|x|$ is well-known to be $\R_+$-homogeneous.

3) The Giseh function $g(x)=d(x,K_\infty)$ where $K_\infty=
\bigcup_{n\in\N}3^n\K$ and $\K$ is the triadic Cantor set, is 
$\frac{1}{3}$-fractal.

\includegraphics{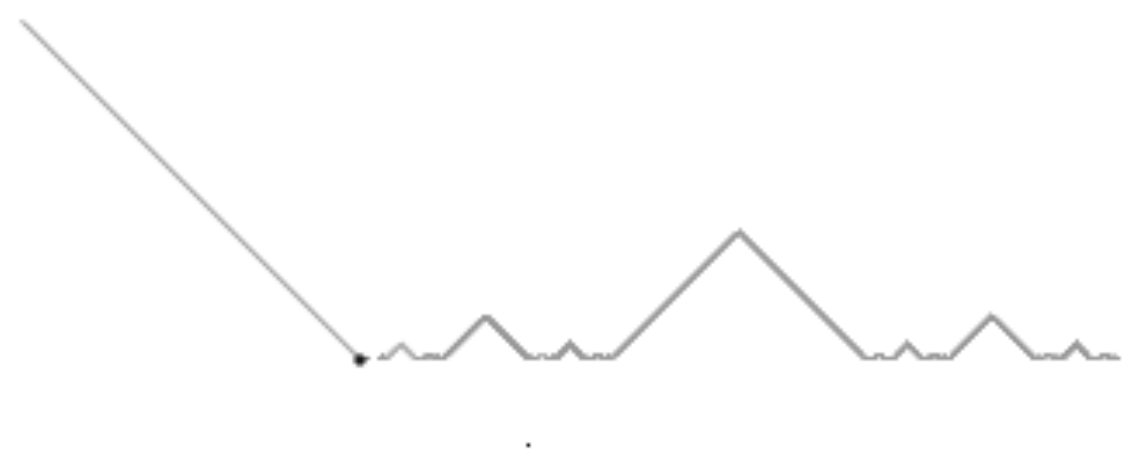}
\begin{picture}(0,5)
\end{picture}

\vspace{-6mm}
3) The fractal wave function  $\xi(x)=x\sin\log|x|$ if $x\not=0$, $\xi(0)=0$, is\break
$e^{-2\pi}$-fractal.

\includegraphics{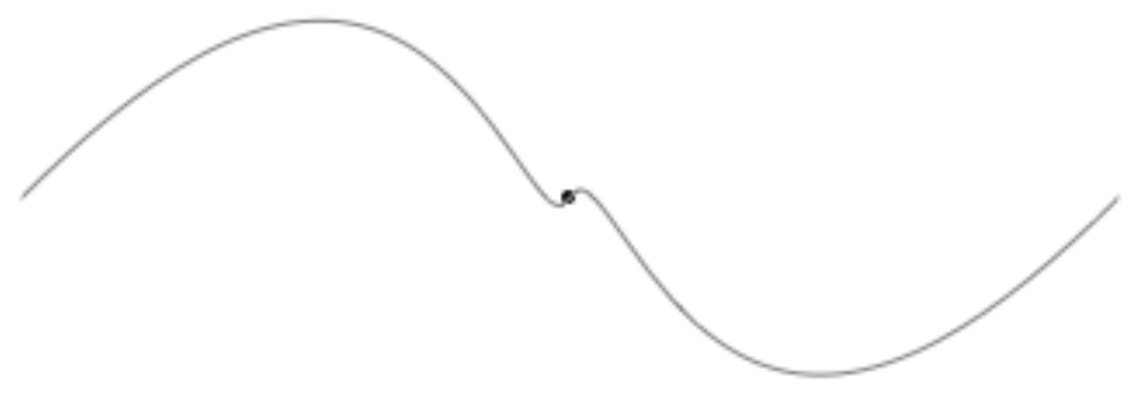}
\begin{picture}(0,0)
\end{picture}

\vspace{3mm}

\begin{prop}
If $h:E\lra E'$ is $\Sigma$-homogeneous, then
$h$ is $v(t_0)$-fractal, where $0<v(t_0)<1$.
\end{prop}

\begin{prop}\rm{($\Sigma$-uniqueness property)}
If $h_1,h_2:E_a\lra E'_{a'}$ are
$\Sigma$-homogeneous; then:
$h_1\tang_a h_2$ $\Lra$
$h_1=h_2$.
\end{prop}

\begin{prop}
If $h:E_a\lra E'_{a'}$ is $\Sigma$-homogeneous, then:\par
\noindent $h$ $LL_a\Llra h$ lipschitzian $\Llra \exists g\ LL_a,\
g\tang_a h\Llra h\ Tang_a$;
and then\break
$\rho(\mathrm{T}h_a)=\inf\{k>0| h\ \, k$-$lipschitzian\}=
\sup_{x\not= y}\frac{\|h(x)-h(y)\|}{\|x-y\|}$. 
\end{prop}

\begin{prop}
If $h:E_a\lra E'_{a'}$ is $\Sigma$-homogeneous, then
$h$ is $\rho(\mathrm{T}h_a)$-lipschitzian, \textit{i.e}
this  lipschitzian ratio $\rho(\mathrm{T}h_a)$ is ``reached'' in $h$.
\end{prop}

$\Sigma$-Lhomogeneous will mean lipschitzian $\Sigma$-homogeneous.
Let us denote\break
$\Sigma$-$\L$hom$(E_a,E'_{a'})=
\{h:E_a\lra E'_{a'}|h$ $\Sigma$-Lhomogeneous, $h(a)=a'\}$;
it is a subset of $\L$L$((E,a),(E',a'))$.
The set $\Sigma$-$\L$hom$(E_a,E'_0)$ is a sub-vector space of
$\L$L$((E,a),(E',0))$. We denote also
$\Sigma$-$\L$hom$(E,E')$ the vector space $\Sigma$-$\L$hom$(E_0,E'_0)$;
of course, $\L(E,E')$ is itself a sub-vector space of $\Sigma$-$\L$hom$(E,E')$ for all $\Sigma$.

Thanks to the $\Sigma$-uniqueness property, 
the restriction of the canonical surjection $q$ to
$\Sigma$-$\L$hom$(E_a,E'_{a'})
\lra\J$et$((E,a),(E',a'))$
is injective,
which allows to define a distance
$d(h,h')=d($T$h_a,$T$h'_a)$ on
$\Sigma$-$\L$hom$(E_a,E'_{a'})$.
We recall that this distance on $\Sigma$-$\L$hom$(E_a,E'_{a'})$ was defined at first as a quasi-distance
on $\L$L$((E,a),(E',a'))$.

\begin{prop}
The above distance on
$\Sigma$-$\L\mathrm{hom}(E_a,E'_{a'})$ 
can be written
$d(h,h')=
\sup_{x\not=a}\frac{\|h(x)-h'(x)\|}{\|x-a\|}$.
\end{prop}

\proof
Denoting $\delta(h,h')=
\sup_{x\not=a}\frac{\|h(x)-h'(x)\|}{\|x-a\|}$ and 
referring to the glossary, we show that, for all $r>0$, we have 
$d^r(h,h')=\delta(h,h')$. Indeed, we have immediately
$d^r(h,h')\leq\delta(h,h')$.
Now, if $x\in E_a$, $t\in \Sigma$ (with $0<v(t)<1$) and $n\in\N$
verify $\|t^n\star_ax\|_a\leq r$ (recalling that $\lim_{n\rightarrow +\infty}
t^n\star_ax=a$),
we use the $\Sigma$-homogeneousness of $h$ to obtain:
$\frac{\|h(x)-h'(x)\|}{\|x-a\|}=\frac{\|h(t^n\star_a x)-h'(t^n\star_a x)\|}{\|t^n\star_a x-a\|}
\leq d^r(h,h')$, which implies $\delta(h,h')\leq d^r(h,h')$.

\begin{prop}
$\Sigma$-$\L\mathrm{hom}(E,E')$ is a n.v.s.
where $\|h\|=
\sup_{x\not=0}\frac{\|h(x)\|}{\|x\|}$.
This norm verifies $\|h(x)\|\leq\|h\|\,\|x\|$
and $\|h'.h\|\leq\|h'\|\,\|h\|$ for all $x\in E$, if\break
$h\in\Sigma$-$\L\mathrm{hom}(E,E')$ and
$h'\in\Sigma$-$\L\mathrm{hom}(E',E'')$.
It goes without saying that\break
$\L(E,E')$, equipped with its operator norm, is a sub-n.v.s. of $\Sigma$-$\L\mathrm{hom}(E,E')$. 
\end{prop}

\proof
We just have to use again the isometric embedding
$\Sigma$-$\L$hom$(E,E')
\buildrel q\over\lra
\J$et$((E,0),(E',0)):h\mapsto$ T$h_0$, which here is also linear, to define the norm\break
$\|h\|=\|$T$h_0\|=d($T$h_0,\calo)=d(h,0)$
on $\Sigma$-$\L$hom$(E,E')$.

\subsection{representable jets}

The lipschitzian homogeneous maps will ``represent'' 
some metric jets, said to be representable; having in mind the example of  
the continuous affine map A$f_a:E\lra E'$ tangent at $a\in U$ to
$f:U\lra U'$ (supposed to be differentiable at $a\in U$;
$U$ and $U'$ being open subsets of $E$ and $E'$ respectively)
which plays the starring role in its tangent jet
T$f_a\in\J$et$((U,a),(U',f(a))$
(or, better said, in its streching jet
$\Gamma($T$f_a$) to $E$): actually, 
A$f_a$ is the unique continuous affine map of the jet
$\Gamma($T$f_a)$; thus, in a way, we could  say
that A$f_a$ ``represents'' the jet $\Gamma($T$f_a)$.

\begin{prop}
The map:
$\Sigma$-$\L\mathrm{hom}(E,E')\lra
\Sigma$-$\L\mathrm{hom}(E_a,E'_{a'}):h\mapsto h^a$, where
$h^a(x)=a'+h(x-a)$, is a translate of $h$ in $a$,
is a bijective isometry.
\end{prop}

A jet $\varphi:(E,a)\lra (E',a')$ is said to be 
$\Sigma$-\textit{representable} if there exists\break 
$h\in\Sigma$-$\L$hom$(E_a,E'_{a'})$ such that $h\in\varphi$
(\textit{i.e}, $h$ being $Tang_a$, T$h_a=\varphi$); 
which is equivalent to say that there exists 
$h\in\Sigma$-$\L$hom$(E,E')$ such that $h^a\in\varphi$.
Such an element of $\varphi$ is unique (thanks to the $\Sigma$-uniqueness property) and is called the\break $\Sigma$-representative element of $\varphi$: it  plays 
a central role in $\varphi$ since, on the one hand $\rho(\varphi)$
is ``reached'' in  it (see prop. 3.6; besides a $\Sigma$-representable jet $\varphi:(E,0)\lra(E',0)$ is good iff
$\rho(\varphi)=\|h\|$) and, on the other hand it gives the ``direction'' of $\varphi$, since
it verifies the following property:

\begin{prop}
If $\varphi:(E,a)\lra (E',a')$ is a $\Sigma$-representable jet and if $h$ is its\break
$\Sigma$-representative element, then, for
all $f\in\varphi$ and all $x\in E$, we have: 

$h(x)=
\lim_{0\not=v(t)\rightarrow 0}t^{-1}\star'_{a'}f(t\star_a x)=
\lim_{0\not=v(t)\rightarrow 0} (f(a)+\frac{f(a+v(t)(x-a))-f(a)}{v(t)})$.
\end{prop}

\begin{theo}
Let $\varphi\in\J\mathrm{et}((E,0),(E',0))$. Then we have the equivalence:\break
$\varphi$ is the jet of a continuous linear map $\Llra$
$\varphi$ is a $\Sigma$-representable linear jet for a $\Sigma$.
Hence, the equivalence:
$\varphi$ is the jet of a continuous linear map $\Llra$
$\varphi$ is a linear jet,\break 
when $\varphi$ is  $\Sigma$-representable for a $\Sigma$. 
\end{theo}

\proof
See the Proof 2 at the end of this paper.

\subsection{Contactable maps}
We now come to the final generalization of differentiable maps:
the $Tang_a$ maps $f$ which are contactable at $a$, their contact at $a$ being a well-chosen map $\kappa f_a$ in the tangential t$f_a$ (this contact is then 
the analogous
of the differential d$f_a$).

In all that follows, we consider a map
$f:U\lra U'$, where $U$ and $U'$ are open subsets respectively 
of $E$ and $E'$, and $a\in U$.
Such an $f$ is said to be\break
$\Sigma$-\textit{contactable} at $a$
(in short $\Sigma$-$Cont_a$), if $f$ is $Tang_a$ and if
the streching jet $\Gamma($T$f_a):(E,a)\lra(E',f(a))$
is a $\Sigma$-representable jet;
which is equivalent to say that there exists 
$h\in\Sigma$-$\L$hom$(E,E')$ such that
$h^a\tang_a f$.
Still thanks to the $\Sigma$-uniqueness property,
these $h$ and $h^a$ are unique,
$h^a$ being the $\Sigma$-representative element of $\Gamma($T$f_a)$, while $h$, denoted $\kappa f_a$, is called
the $\Sigma$-\textit{contact} of $f$
at $a$; this $\kappa f_a$ is the $\Sigma$-representative element of the tangential of $f$ at $a$, 
t$f_a=\Omega(\Gamma($T$f_a))$.
Referring to prop. 3.10,
this $\Sigma$-contact can be written
$\kappa f_a(x)=
\lim_{0\not=v(t)\rightarrow 0} (\frac{f(a+v(t)x)-f(a)}{v(t)})$
for all $x\in E$.

We still have a composition of $\Sigma$-contacts for composable $\Sigma$-contactable maps:
$\kappa(g.f)_a=\kappa g_{f(a)}.\kappa f_a$ if $f$ is 
$\Sigma$-$Cont_a$ and $g$ is $\Sigma$-$Cont_{f(a)}$.

We have the implication:
$f$ $Dif\! f_a\Lra f$ $\ \Sigma$-$Cont_a$ for all $\Sigma$
with $\kappa f_a=$ d$f_a$
(the inverse implication being true iff the $\Sigma$-contact $\kappa f_a$ is linear).
In fact, applying thm 3.11 to the tangent jet T$(\kappa f_a)_0$,
we have the following result:

\begin{theo}
We have the equivalence:
$f$ $Dif\! f_a$ $\Llra \ f$ $TL_a$ and $f$ $\Sigma$-$Cont_a$ for a $\Sigma$. Hence the equivalence:
$f$ $Dif\! f_a$ $\Llra f\ TL_a$, when $f$ is $\Sigma$-$Cont_a$ for a $\Sigma$.
\end{theo}

We denote $\kappa_+ f_a$
the $\R_+$-contact at $a$ 
of a map $f$ which is $\R_+$-contactable at $a$ (this $\R_+$-contact
is a $\R_+$-Lhomogeneous map).

We say that a map $f$ is $k$-\textit{neofractal} at $a$ 
(in short $k$-neofract$_a$) if it is\break
$\N'_k$-$Cont_a$; its
$\N'_k$-contact at $a$ (denoted $\kappa_k f_a$) is a $k$-\textit{Lfractal} map 
(\textit {i.e} a\break
$\N'_k$-Lhomogeneous map); we denote
$k$-$\L$fract$(E,E')$ the n.v.s. $\N'_k$-$\L$hom$(E,E')$,
neofract$_a$ (resp. Lfractal and $\L$fract$(E,E')$) meaning that such a 
$k\in]0,1[$ exists.

\begin{remks}\par\hfill

1) Of course, by definition of the contactibility, we have the implication:\par
\noindent $h:E_0\lra E'_0$ $\ \Sigma$-Lhomogeneous $\ \Lra\ $ 
$h$ $\ \Sigma$-$Cont_0$, 
with $\kappa h_0=h$. 
In particular,\break
$h$ $\ k$-Lfractal $\Lra$ 
$h$ $\ k$-neofract$_0$, with $\kappa_k h_0=h$.

2) Referring to prop. 3.3, we have the implication:\par
\noindent $f$ $\Sigma$-$Cont_a$ $\Lra f\ $ neofract$_a$.

3) Referring to previous implications, we have the particular implications:\par
\noindent $f$ $Dif\! f_a\Lra f \ \R_+$-$Cont_a\Lra f\ $ neofract$_a$ 
\end {remks}

\begin{exams}
{}
\end{exams}

1) $E$ being a given n.v.s., every norm on $E$ (which is 
equivalent to the given norm on $E$) is $\R_+$-$Cont_0$ with
$\kappa_+ n_0=n$, since $n$ is $\R_+$-Lhomogeneous;
it is the case for the function $v(x)=|x|$.

2) The Giseh function $g(x)=d(x,K_\infty)$, where $K_\infty=
\bigcup_{n\in\N}3^n\K$ and $\K$ is the triadic Cantor set,
is $\frac{1}{3}$-neofract$_0$ 
(even $\frac{1}{3}$-Lfractal: see example 3.2) with $\kappa_{\frac{1}{3}}g_0=g$.

3) The fractal wave function $\xi(x)=x\sin\log|x|$ if $x\not=0$, $\xi(0)=0$, is\break
$e^{-2\pi}$-neofract$_0$ (even
$e^{-2\pi}$-Lfractal: see examples 3.2) with
$\kappa_{e^{-2\pi}}\xi_0=\xi$.
However, this function is not $\R_+$-$Cont_0$, since 
$\kappa_{e^{-2\pi}}\xi_0=\xi$ is not
$\R_+$-Lhomogeneous.

4) Let us consider the following bifractal wave function
defined by\break
$\zeta(x)=x\sin\frac{2\pi}{a}\log|x|$ if $x<0$,
$\zeta(x)=x\sin\frac{2\pi}{b}\log|x|$ if $x>0$, $\zeta(0)=0$, where $a,b$ are $>0$ real numbers verifying $\frac{a}{b}\not\in\Q$.
We will prove (in Proof 3 at the end of this paper)
that this function is $Tang_0$, but not neofract$_0$
(although for $r>0$, $\zeta_r(x)=x\sin\frac{2\pi}{r}\log|x|$ if $x\not=0$, $\zeta_r(0)=0$, is $e^{-r}$-Lfractal!);
thus, referring to remarks 3.13, this 
bifractal wave function is $\Sigma$-$Cont_0$ for none $\Sigma$.

5) Let us at last notice that our uncanny function
$Ins(x)=x\sin\log|\log|x||$ if $x\not=0$, $Ins(0)=0$,
is also not $\Sigma$-$Cont_0$ for any $\Sigma$, since it is
$TL_0$ (see examples 2.4 and examples 2.5) and not $Dif\! f_0$: it thus remains to use thm. 3.12!

\vspace{3mm}
Finaly, the fractal waves allows us to establish the 
following remarkable result
(remarkable, since it deals with jets at order 1!)):

\begin{theo}

$\J\mathrm{et}((\R,0),(\R,0))$ is a n.v.s of infinite dimension.
\end{theo}

\proof
Indeed, by a constructing procedure which is analogous to 
the one of our $e^{-2\pi}$- fractal wave $\xi(x)=x\sin\log|x|$
if $x\not=0$, $\xi(0)=0$, we associate a $e^{-T}$-fractal wave
$\tilde f(x)=xf(\log|x|)$, $\tilde f(0)=0$, to every function 
$f:\R\lra\R$\break
which is lipschitzian, $T$ periodic ($T>0$) and for which there exists a right derivative at each point. Denoting 
$\calp er(T)$ the vector space of these functions $f$, we have an evident embedding $\calp er(T)\lra 
e^{-T}$-$\L$fract$(\R,\R):f\mapsto\tilde f$. 
We just\break
have now to compose this embedding 
with our well-known embedding\break
$q:e^{-T}$-$\L$fract$(\R,\R)\lra\J$et$((\R,0),(\R,0)):
h\mapsto$ T$h_0$ to
obtain an embedding
$\calp er(T)\lra\J$et$((\R,0),(\R,0))$; it remains then to use the fact that
$\calp er(T)$ is of infinite dimension.

\section{Contactibility with some classical Theorems}

Skimming through the metric Dif\!\! ferential Calculus
has highlighted many generalizations of the specific properties of the 
classical differentials.

Actually, for contactable maps, we
can add to these generalizations a mean value theorem;
and theorems about extrema which, unlike the classical ones,
need hypothesis only at order 1!

\vspace{3mm}
In what follows, $U$  and $U'$ are open subsets of 
n.v.s. $E$ and $E'$ respectively.

\begin{theo}
Let $f:U\lra U'$ be
a continuous map, $a,b\in U$ such that $[a,b]\subset U$, $F$ a finite subset of $]a,b[$.
We assume that, for all $x\in\,]a,b[-F$, the map $f$ is\break
$\Sigma$-$Cont_x$  and satisfies 
$\|\kappa f_x\|\leq k$ (where $k\geq 0$; the previous norm has been defined in prop.3.8). Then
we have: $\|f(b)-f(a)\|\leq k\|b-a\|$.
\end{theo}

\begin{theo}
Let $f:U\lra\R$ be
$\Sigma$-$Cont_a$ (with $a\in U$), admitting a local minimum at $a$.
Then $\kappa f_a$ admits a global minimum at 0.
\end{theo}

\vfill\eject

\begin{remk}

This gives back the well-known result of the differentiable
case:\break
`` $f$ admits a local minimum at $a\Lra a$ is a critical point of $f$ (i.e $\textup{d}f_a=0$)''.
Indeed, if $f$ is $Dif\! f_a$, then $f$ is $Cont_a$ with 
$\kappa f_a=\mathrm{d}f_a$, so that the $\Sigma$-contact 
$\kappa f_a$ is a continuous linear function $E\lra\R$
admitting a global minimum at 0: it forces this $\kappa f_a$ to be the null function. 
\end{remk}

\begin{theo}
Let $f:U\lra\R$ be
$\R_+$-$Cont_a$ (with $a\in U$; $E$ being here of finite dimension),
such that $\kappa_+ f_a>0$ (i.e verifying $\kappa_+ f_a(x)>0$ for every $x\in E-\{a\}$).
Then, $f$ admits a strict local minimum at $a$.
\end{theo}

\begin{remk}
This theorem has not its equivalent, at order 1, in 
classical Differential Calculus, since a linear function cannot have a strict minimum. It is rather inspired by theorems giving sufficient conditions, at order 2, 
for the existence of extrema.
\end{remk}

\vspace{10mm}
\centerline{\textbf{PROOFS}}
\vspace{6mm}

{\bf Proof 1}:
\textit{
We prove here that the jet $\cali$ of the uncanny function
$Ins$, defined in examples 2.1
is a good jet,
with $\rho(\cali)=1$.
}

First, we notice that $Ins$ is 2-$LL_0$ (since $Ins$ is odd and, on
$]0,\frac{1}{e}[$, it verifies $|Ins'(x)|\leq
1+\frac{1}{|\log x|}\leq 2$. Thus $\rho(\cali)\leq 2$;
and even $\rho(\cali)\leq 1+\frac{1}{\alpha}$ for all
$\alpha>0$, since $Ins$ is in fact $(1+\frac{1}{\alpha})$-lipschitzian on $]-e^{-\alpha},e^{-\alpha}[$). 
Finally, we have $\rho(\cali)\leq 1$, and thus
$d(\cali,\calo)\leq\rho(\cali)\leq 1$.
It remains to prove that $1\leq d(\cali,\calo)$ \textit{i.e}
that $1\leq d(Ins,0)$ (since $\cali=$T$Ins_0$ and
$\calo=$T$0_0$); for this, we use the definition of the quasi-distance recalled in the below glossary.
Well, by definition, $d(Ins,0)=
\lim_{x\rightarrow 0}d^r(Ins,0)$,
where $d^r(Ins,0)=\sup\{\frac{|Ins(x)|}{|x|}\ |\ 0<|x|\leq r\}$.
We notice that, for all $r>0$, there exists $k$
verifying $x_k=e^{-e^{\frac{\pi}{2}+2k\pi}}\leq r$,
so that $|\frac{Ins(x_k)}{x_k}|=|\sin\log|\log x_k||=1$
which implies $d^r(Ins,0)\geq 1$ for all $r>0$.
It remains then to do $r\rightarrow 0$.

\vspace{3mm}
{\bf Proof 2}:
\textit{We prove here the theorem 3.11.}

Let $h$ be the $\Sigma$-representative element of a 
$\Sigma$-representable jet $\varphi:\break
(E,0)\lra(E',0)$; thus $\varphi=$T$h_0$
with $h\in\L$hom$(E,E')$.
Let us also assume that this jet is $TL_0$; we can then write
T$(h.\sigma)_{(0,0)}=$T$h_0.$T$\sigma_{(0,0)}=
\varphi.+=+.\varphi^2=$
T$\sigma_{(0,0)}.$T$h_0^2=$T$(\sigma.h^2)_{(0,0)}$,
\textit{i.e}
$h.\sigma\tang_{(0,0)}\sigma.h^2$.
Now, by the $\Sigma$-uniqueness property, we deduce
the equality $h.\sigma=\sigma.h^2$, 
which gives the linearity of $h$ (since $h$ is continuous);
that is the end of this proof, since $h\in\varphi$
(the inverse implication being evident).

\vspace{3mm}
{\bf Proof 3}:
\textit{We prove here that the bifractal wave function
defined in examples 3.14, is $Tang_0$ but not neofract$_0$.}

Let us consider first the function  
$\zeta_r(x)=x\sin\frac{2\pi}{r}\log|x|$ if $x\not=0$, $\zeta_r(0)=0$ (where $r>0$); it is derivable on $\R^*$
with $|\zeta'_r(x)|\leq k_r=1+\frac{2\pi}{r}$;
thus our bifractal wave function $\zeta$ is $\sup(k_a,k_b)$-lipschitzian which implies that it is 
$Tang_0$.

We prove now that $\zeta$ cannot be neofract$_0$.
Indeed, if there exists $k\in ]0,1[$ for which $\zeta$ is 
$k$-neofract$_0$, then there exists a $k$-Lfractal function $g:\R\lra\R$\break
verifying $g\tang_0 \zeta$.
Thus, for all $x\in \R^*$, we have (since 
$\lim_{n\rightarrow  \infty}k^nx=0$):
$\lim_{n\rightarrow\infty}\frac{|g(k^n x)-\zeta(k^n x)|}
{|k^n x|}=0$. Using the fact that $g$ is $k$-fractal, we can write
$g(k^n x)=k^ng(x)$ for all $n\in\N$, so that
$\lim_{n\rightarrow\infty}|\frac{g(x)}{x}-
\frac{\zeta(k^n x)}{k^n x}|=0$ which gives
$\lim_{n\rightarrow\infty}\frac{\zeta(k^n x)}{k^n x}=
\frac{g(x)}{x}$.
Now, for $x<0$, $\frac{\zeta(k^n x)}{k^n x}=
\sin(n\frac{2\pi\log k}{a}+\frac{2\pi}{a}\log|x|)$;
so, if we set $\alpha=\frac{\log k}{a}$, 
$\gamma=\frac{2\pi}{a}\log|x|$ and $x_n=\sin(2\pi n\alpha+\gamma)$,
we have obtained that the sequence $(x_n)$ converges
towards $\frac{g(x)}{x}$.
So the cluster set of the sequence $(x_n)$ is reduced to $\frac{g(x)}{x}$
and thus cannot be equal to $[-1,1]$; this implies that $\alpha\in\Q$
(see the lemma 4.6 below). 
In the same way (when $x>0$), we show that $\beta=\frac{\log k}{b}\in \Q$. Then $\frac{a}{b}=\frac{\beta}{\alpha}\in\Q$,
which contredicts the hypothesis made on the definition of $\zeta$! Thus, such a $k$ cannot exist.

\begin{lema}

Let $\alpha\in\Q^c$ and $\gamma\in\R$. We recall the following results:

1) The set $\{e^{2\pi in\alpha}\, |\, n\in\Z\}$ is dense in 
$S^1=\{z\in\C\, |\, |z|=1\}$; and even,  the set $\{e^{2\pi in\alpha}\, |\, n\in\N\}$ is still dense in $S^1$.

2) Using the bijective isometry
$\varphi:\C\lra\C:z\mapsto e^{i\gamma}z$,
we obtain that the set $\{e^{i(2\pi n\alpha+\gamma)}\, |\, n\in\N\}$ is 
also dense in $S^1$.

3) Setting $x_n=\sin(2\pi n\alpha+\gamma)$ for all $n\in\N$, we deduce, not only that the set $\{x_n\, |\, n\in\N\}$ is dense in $[-1,1]$, but, even better, that the cluster set
of the sequence $(x_n)_{n\in\N}$ is equal to $[-1,1]$.

\end{lema}
\vspace{6mm}
\centerline{\textbf{GLOSSARY}}
\vspace{4mm}

n.v.s.: $\R$-normed vector space, usually denoted $E$.

\vspace{1mm}

p.m.s.: pointed metric space, usually denoted $(M,a)$ (where 
$a$ is not isolated).

\vspace{1mm}
$C^0_a$: continuous at $a$ ($C^0$: continuous).

\vspace{1mm}

$LL_a$: locally lipschitzian at $a$.

\vspace{1mm}
$LSL_a$: locally semi-lipschitzian at $a$ (knowing that $f:M\lra M'$ is $SL_a$ if\break
$\exists k>0\ \forall x\in M\ d(f(x),f(a))\leq kd(x,a)$).
We have the implications:\break
$f\ LL_a\Lra f\ LSL_a\Lra f\ C^0_a$.

\vspace{1mm}
$\L$L$((M,a),(M',a'))$: the set of maps $f:M\lra M'$ which are $LL_a$ and which verify $f(a)=a'$. These sets are the ``Hom''
of a cartesian category $\L$L. 

\vspace{1mm}

$\J$et$((M,a),(M',a'))=\L$L$((M,a),(M',a'))/\tang_a$;
these sets are the ``Hom''
of the cartesian category $\J$et (whose objects are the p.m.s. and the morphims the metric jets (or jets)). This category $\J$et is a quotient of the category $\L$L by the relation of tangency. The canonical surjection $q:\L{\rm L}\lra\J{\rm{et}}$ is a cartesian functor.

\vspace{1mm}

$\M$et: the category whose objects are the metric spaces and morphisms
the $LSL$ maps. $\J$et is enriched in $\M$et:
for $f,g\in\L$L$((M,a),(M',a'))$, there exists a $k>0$ and a neighborhood $V$ of $a$ on which $f$ and $g$ are $k$-lipschitzian ; we first define
$d(f,g)=\lim_{r\rightarrow 0}d^r(f,g)$, where $d^r(f,g)=\sup
\{\frac{d(f(x),g(x))}{d(x,a)}\ |\ a\not=x\in B'(a,r)\cap V\}$.
It is not a distance since we only have $d(f,g)=0\Llra f\tang_a g$ (it is only a ``quasi-distance''); this leads us to define a true distance
on the quotient $\J$et$((M,a),(M',a'))$ by setting
$d(q(f),q(g))=d(f,g)$.

\vspace{1mm}
$Tang_a$: tangentiable at $a$; T$f_a$ being the tangent jet at $a$ of $f$ (assumed to be $Tang_a$).
If $f$ is $LL_a$, we have $f\in$T$f_a=q(f)$.

\vspace{6mm}
\noindent In what follows, $f:U\lra U'$, where $U$ and $U'$ are open subsets of n.v.s. $E$ and
$E'$ respectively; and $a\in U$. 

\vspace{2mm}
$Dif\! f_a$: differentiable at $a$; d$f_a$ being the differential at $a$
of $f$ (assumed to be $Dif\! f_a$).

\vspace{1mm}
t$f_a$: tangential at $a$ of $f$ (assumed to be $Tang_a$).

\vspace{1mm}
$CT$: continuously tangentiable, \textit{i.e} tangentiable at every point of $U$ and the map
t$f:U\lra\J{\rm{et}}((E,0),(E',0)):x\mapsto$ t$f_x$, called the tangential of $f$, 
is continuous.

\vspace{1mm}
$TL_a$: tangentially linear at $a$.

\vspace{1mm}
$CTL$: continuously tangentially linear \textit{i.e} tangentially linear at every point of $U$
and the restriction
$\lambda f$ of t$f$ to $U\lra\Lambda(E,E')$ is continuous.

\noindent We have the implications:
$\quad C^1\Lra CTL\Lra CT\Lra C^0$.

\vspace{1mm}
$\L(E,E')$: the set of all continuous linear maps $E\lra E'$.

\vspace{1mm}
$\Lambda(E,E')$: the set of all linear jets $(E,0)\lra (E',0)$.

\vspace{1mm}
$\C$TL: the category whose objects are the open subsets of n.v.s.
and whose morphisms are the $CTL$ maps.

\vspace{1mm}
$\G$CTL: the category whose objects are the n.v.s., and morphisms $E\lra E'$ are germs at 0, of maps $f:E\lra E'$ verifying  $f(0)=0$ and for which there exists a neigborhood $V$ of 0 such that $f|_V:V\lra E'$ is $CTL$.


\vspace{1mm}
$\Sigma$: valued monoid (its valuation being denoted $v:\Sigma\lra\R_+$).

\vspace{1mm}
$\Sigma$-$Cont_a$: $\Sigma$-contactable at $a$.

\vspace{3mm}
$\Sigma$-contracting space: metric space $M$, centred in $\omega$, on which $\Sigma$ externally operates, this operation 
verifying $0\star x=\omega$ for all $x\in  M$ and 
$t\star \omega=\omega$ for all $t\in\Sigma$;
and being compatible with the distance of $M$,
\textit{i.e} verifying
$d(t\star x,t\star y)=v(t)d(x,y)$ for all
$(t,x,y)\in\Sigma\times M\times M$.


\end{document}

%% file: macro.tex
\def \lra{\longrightarrow}

\def\Lra{\Longrightarrow}

\def\Llra{\Longleftrightarrow}

\def\proof{\vspace{0mm}$\underline{\hbox{\sl Proof}}$ :\ }

\newtheorem{thm}{Théorème}[section]

\newtheorem{theo}[thm]{Theorem}
\newtheorem{prop}[thm]{Proposition}

\newtheorem{lema}[thm]{Lemma}

\newtheorem{remk}[thm]{Remark}
\newtheorem{remks}[thm]{Remarks}

\newtheorem{exams}[thm]{Examples}

\def\R{\mathbb R}
\def\L{\mathbb L}

\def\J{\mathbb J}
\def\M{\mathbb M}
\def\C{\mathbb C}

\def\N{\mathbb N}
\def\Z{\mathbb Z}

\def\K{\mathbb K}
\def\Q{\mathbb Q}
\def\G{\mathbb G}

\def\calf{{\mathcal F}}
\def\calg{{\mathcal G}}

\def\cali{{\mathcal I}}
\def\calj{{\mathcal J}}

\def\calo{{\mathcal O}}
\def\calp{{\mathcal P}}

\def\calv{{\mathcal V}}
